\documentclass[12pt]{article}

\usepackage[latin1]{inputenc}
\usepackage{bbm}
\usepackage{stmaryrd}
\usepackage{latexsym}
\usepackage{amsfonts}
\usepackage{amsmath,amssymb,amscd}
\usepackage{yfonts}
\usepackage{dsfont}
\usepackage{mathabx}
\usepackage[dvips]{graphicx}
\usepackage{epsfig}
\usepackage{psfrag}
\usepackage[font={small}]{caption}

\DeclareFontFamily{U}{txsyc}{}
\DeclareFontShape{U}{txsyc}{m}{n}{
   <-> txsyc%
}{}
\DeclareFontShape{U}{txsyc}{bx}{n}{
   <-> txbsyc%
}{}
\DeclareFontShape{U}{txsyc}{l}{n}{<->ssub * txsyc/m/n}{}
\DeclareFontShape{U}{txsyc}{b}{n}{<->ssub * txsyc/bx/n}{}
\DeclareSymbolFont{symbolsC}{U}{txsyc}{m}{n}
\SetSymbolFont{symbolsC}{bold}{U}{txsyc}{bx}{n}
\DeclareFontSubstitution{U}{txsyc}{m}{n}
\DeclareMathSymbol{\df}{\mathrel}{symbolsC}{"42}
\DeclareMathSymbol{\fd}{\mathrel}{symbolsC}{"43}
\DeclareMathSymbol{\lJoin}{\mathrel}{symbolsC}{"58}
\DeclareMathSymbol{\rJoin}{\mathrel}{symbolsC}{"59}

\newcommand{\cE}{{\cal E}}
\newcommand{\cF}{{\cal F}}

\newcommand{\cL}{{\cal L}}

\newcommand{\cO}{{\cal O}}
\newcommand{\cP}{{\cal P}}

\newcommand{\EE}{\mathbb{E}}

\newcommand{\LL}{\mathbb{L}}
\newcommand{\NN}{\mathbb{N}}
\newcommand{\PP}{\mathbb{P}}

\newcommand{\RR}{\mathbb{R}}

\newcommand{\ZZ}{\mathbb{Z}}

\newcommand{\iy}{\infty}
\newcommand{\lt}{\left}
\newcommand{\me}{\medskip}

\newcommand{\pa}{\partial}
\newcommand{\ri}{\rightarrow}
\newcommand{\rt}{\right}
\newcommand{\sm}{\smallskip}

\newcommand{\wi}{\widetilde}
\newcommand{\wit}{\widehat}

\newcommand{\ex}{\exists\ }
\newcommand{\fo}{\forall\ }

\newcommand{\lve}{\lt\vert}
\newcommand{\lVe}{\lt\Vert}

\newcommand{\rve}{\rt\vert}
\newcommand{\rVe}{\rt\Vert}

\newcommand{\st}{\,:\,}

\newcommand{\un}{\mathds{1}}

\newcommand{\bq}{\begin{eqnarray*}}
\newcommand{\bqn}[1]{\begin{eqnarray}\label{#1}}
\newcommand{\eq}{\end{eqnarray*}}
\newcommand{\eqn}{\end{eqnarray}}
\newcommand{\wwtbp}{\par\hfill $\blacksquare$\par\me\noindent}
\newcommand{\thistitlepagestyle}{}
\newcommand{\lin}{\llbracket}
\newcommand{\rin}{\rrbracket}

\newcommand{\ttsim}{\raise.17ex\hbox{$\scriptstyle\mathtt{\sim}$}}

\newtheorem{pro}{Proposition} 

\newtheorem{lem}[pro]{Lemma}
\newtheorem{theo}[pro]{Theorem}
\renewcommand{\thepro}{\arabic{pro}}

\newenvironment{exa}
{\par\me\refstepcounter{pro}\noindent{\bf Example \thepro\ }}
{\par\hfill $\square$\par\sm\noindent}
\newenvironment{exas}
{\par\me\refstepcounter{pro}\noindent{\bf Example \thepro\ }}
{\par\hfill $\square$\par\sm\noindent}

\newcommand{\proof}{\par\me\noindent\textbf{Proof}\par\sm\noindent}
\newcommand{\prooff}[1]{\par\me\noindent\textbf{#1}\par\sm\noindent}

\newcommand{\comment}[1]{}

\title{Estimates on the amplitude of the first Dirichlet eigenvector\\ in discrete frameworks}
\author{Persi Diaconis\footnote{Partially supported by ANR-11-LABX-0040-CIMI within the program ANR-11-IDEX-0002-02} \ and Laurent Miclo\footnote{Partially supported by ANR-12-BS01-0019}
}
\date{\box1
 \box2
}

\begin{document}

\setbox1=\vbox{
\large
\begin{center}
${}^\dagger$Department of Mathematics\\
Stanford University, California, USA\\
\end{center}
} 
\setbox2=\vbox{
\large
\begin{center}
${}^\ddagger$Institut de Mathématiques de Toulouse, UMR 5219\\
Universit\'e de Toulouse and CNRS, France\\
\end{center}
} 
\setbox3=\vbox{
\hbox{Stanford University\\}
\hbox{Department of Mathematics\\}
\hbox{Building 380, Sloan Hall\\}
\hbox{Stanford, California 94305, USA\\}
}
\setbox4=\vbox{
\hbox{${}^\ddagger$miclo@math.univ-toulouse.fr\\}
\vskip1mm
\hbox{Institut de Mathématiques de Toulouse\\}
\hbox{Université Paul Sabatier\\}
\hbox{118, route de Narbonne\\} 
\hbox{31062 Toulouse Cedex 9, France\\}
}
\setbox5=\vbox{
\box3
\vskip5mm
\box4
}

\maketitle
\thistitlepagestyle
\abstract{
Consider a finite absorbing Markov generator, irreducible on the non-absorbing states. 
Perron-Frobenius theory ensures the existence of a corresponding positive eigenvector $\varphi$.
The goal of the paper is to give bounds on the amplitude $\max \varphi/\min\varphi$.
Two approaches are proposed: one using a path method and the other one, restricted to the reversible situation,
based on spectral estimates. The latter approach is extended to denumerable birth and death processes absorbing at 0  for
which infinity is an entrance boundary.
The interest of estimating the ratio is the reduction of the quantitative study of convergence
to quasi-stationarity to the convergence to equilibrium of related ergodic processes, as  seen in \cite{diaconis:hal-01002622}.
}
\vfill\null
{\small
}\par

\newpage 

\section{Introduction}

This paper, a companion to \cite{diaconis:hal-01002622}, develops tools to get useful quantitative bounds on rates of convergence to quasi-stationarity 
for absorbing Markov processes. 
With notation explained below, the bounds in \cite{diaconis:hal-01002622} are of the form
\bq
 \frac{\varphi_\wedge}{2\varphi_\vee}\lVe \wi\mu_0\wi P_t -\wi\eta\rVe_{\mathrm{tv}}\ \leq \ \lVe \mu_t-\nu\rVe_{\mathrm{tv}}\ \leq \ 2\frac{\varphi_\vee}{\varphi_\wedge}\lVe \wi\mu_0\wi P_t -\wi\eta\rVe_{\mathrm{tv}}.\eq
In the middle is the term of interest: $\mu_t$ is the transition probability conditioned on non-absorbtion  at time $t\geq 0$ and $\nu$ is the quasi-stationary distribution.
On both sides, $\wi P_t$ is the Doob transform (forced to be non-absorbing), $\wi\mu_0$ is an associated starting distribution and $\wi\eta$ is the stationary distribution of the transformed process. The point is that quantitative rates of convergence to quasi-stationarity are hard to come by, requiring new tools which are not readily available. The pair
$(\wi P_t,\wi\eta)$ is a usual ergodic Markov chain with many techniques available.\\
The two sides differ by a factor ${\varphi_\wedge}/{2\varphi_\vee}$. Here $\varphi$ is the usual Peron-Forbenius eigenfunction for the matrix restricted to the non-absorbing sites and $\varphi_\wedge\df \min \varphi$, $\varphi_\vee\df \max \varphi$. For the bounds to be useful, we must get control of this ratio. In  \cite{diaconis:hal-01002622}, this control was achieved in special examples where analytic expressions are available with explicit diagonalization. The purpose of the present paper is to give a probabilistic interpretation of this ratio as well as several bounding techniques. For background on quasi-stationarity see M{\'e}l{\'e}ard and Villemonais \cite{MR2994898}, Collet, Mart{\'{\i}}nez and  San~Mart{\'{\i}}n \cite{MR2986807}, van Doorn and  Pollett \cite{MR3063313}, Champagnat and Villemonais \cite{2014arXiv1404.1349C}   or the discussion in \cite{diaconis:hal-01002622}.
We proceed to a more careful description.\par\me
Let us begin by introducing the finite setting. 
 The whole finite state space is $\bar S\df S\sqcup\{\iy\}$, where $\iy$ is the absorbing point.
 This means that $\bar S$ is endowed with a Markov generator matrix $\bar L\df(\bar L(x,y))_{x,y\in\bar S}$
 whose restriction to $S\times S$ is irreducible and  such that
 \bq
 \fo x\in \bar S,\qquad \bar L(\iy,x)&=&0\\
 \ex x\in S \st\qquad \bar L(x,\iy)&>&0.\eq
 Recall that a Markov (respectively subMarkovian) generator is a matrix whose off-diagonal entries are non-negative and 
 such that the sums of the entries of a row all vanish (resp.\ are non-positive).
 \par
 An eigenvalue $\lambda$ of $\bar L$ is said to be of Dirichlet type if an associated eigenvector vanishes at $\iy$.
 Equivalently, $\lambda$ is an eigenvalue of the $S\times S$ minor $K$ of $\bar L$.
Since the matrix $K$
is an irreducible  subMarkovian  generator, the Perron-Frobenius theorem implies
that $K$ admits a unique eigenvalue $\lambda_0$ whose associated eigenvector is positive.
The eigenvalue $\lambda_0$ is simple and we denote by $\varphi$ an associated positive eigenvector.
Its renormalization is not very important for us, because we will be mainly concerned by its \textit{amplitude}
defined by
\bq
a_\varphi&\df&\frac{\varphi_\vee}{\varphi_\wedge},
\eq
with
\bq
\varphi_\vee\ \df\ \max_{x\in S}\varphi(x),\qquad \varphi_\wedge\ \df\ \min_{x\in S}\varphi(x).\eq
\par
We refer to \cite{diaconis:hal-01002622} for the importance of $a_\varphi$ in the investigation of the convergence to quasi-stationarity of
the absorbing Markov processes generated by $\bar L$. Our purpose here is to
estimate this quantity.\par\me
Our approach is based on a probabilistic interpretation of $\varphi$ and, more precisely,
of the ratios of its values.
For any $x\in S$, let $X^x\df(X^x_t)_{t\geq 0}$ be a càdlàg Markov process generated by $\bar L$ and starting from $x$.
For any $y\in \bar S$, denote by $\tau_y^x$ the first hitting time of $y$ by $X^x$:
\bqn{tauxy}
\tau^x_y&\df& \inf\{t\geq 0\st X^x_t=y\},\eqn
with the convention that $\tau_y^x=+\iy$ if $X^x$ never reaches $y$. The first identity below comes from Jacka and Roberts \cite{MR1363332}.
\begin{pro}\label{pro0}
For any $x,y\in S$, we have
\bq
\frac{\varphi(x)}{\varphi(y)}
&=&\EE[\exp(\lambda_0\tau_y^x)\un_{\tau^x_y<\tau^x_\iy}].\eq
In particular, with $O\df\{x\in S\st \bar L(x,\iy)>0\}$,  we have
\bq
a_\varphi&=&\max_{x\in S,\, y\in O}\EE[\exp(\lambda_0\tau_y^x)\un_{\tau^x_y<\tau^x_\iy}].\eq
\end{pro}\par
This probabilistic interpretation leads to two methods of estimating $a_\varphi$.
The first one is through a path argument.
\par
If $\gamma=(\gamma_0,\gamma_1, ..., \gamma_l)$ is a path in $S$, with $\bar L(\gamma_k,\gamma_{k+1})>0$ for all $k\in\lin 0, l-1\rin$, denote
\bqn{Pg}
P(\gamma)&\df& \prod_{k\in \lin 0, l-1\rin}\frac{\bar L(\gamma_k,\gamma_{k+1})}{\lve \bar L(\gamma_k,\gamma_k)\rve-\lambda_0}
\eqn
(for any $l'\leq l''\in\ZZ$, $\lin l',l''\rin\df\{l', l'+1, ..., l''-1, l''\}$ and for $l'\in \NN$,  $\lin l'\rin\df \lin 1,l'\rin$).
\begin{pro}\label{pro1}
Assume that for any $y\in O$ and $x\in S$, we are given a path $\gamma_{y,x}$ going from $y$ to $x$.
Then we have
\bq
a_\varphi&\leq & \lt(\min_{y\in O,\, x\in S} P(\gamma_{y,x})\rt)^{-1}.\eq
\end{pro}
\par
The second method requires that $K$ (the generator restricted to $S$) admit a reversible probability $\eta$ on $S$, namely satisfying
\bq
\fo x,y\in S,\qquad \eta(x)\bar L(x,y)&=&\eta(y)\bar L(y,x).\eq
The operator $-K$ is then diagonalizable. Let $\lambda_0<\lambda_1\leq \lambda_2\leq \cdots\leq \lambda_{N-1}$ be
its eigenvalues, where $N$ is the cardinality of $S$ (the first inequality is strict, due to the Perron Frobenius theorem
and to the irreducibility of $K$).
For any $x\in S$, let $\lambda_0(S\setminus\{x\})$ be the first eigenvalue of the $(S\setminus\{x\})\times (S\setminus\{x\})$
minor of $-K$ (or of $-\bar L$). Finally, consider
\bqn{wila}
 \lambda_0'&\df& \min_{x\in O} \lambda_0(S\setminus\{x\}).\eqn
\begin{pro}\label{pro2}
Under the reversibility assumption, we have 
\bq
a_\varphi&\leq &\lt(\lt( 1-\frac{\lambda_0}{\lambda'_0}\rt)
\prod_{k\in \lin N-1\rin} \lt( 1-\frac{\lambda_0}{\lambda_k}\rt)\rt)^{-1}.\eq
\end{pro}
\par
One advantage of the last result is that it can be extended to absorbing processes on denumerable state spaces,
at least under appropriate assumptions.
We won't develop a whole theory here, so let us just give the example
 of birth and death processes on $\ZZ_+$ absorbing at 0 and for which $\iy$ is an entrance boundary.
 To follow the usual terminology in this domain, we  change the notation, 0 being the absorbing point
 and $\iy$ being the boundary point at infinity of $\ZZ_+$.
We consider $S\df\NN\df\{1,2,3, ...\}$ and $\bar S\df\ZZ_+\df \{0,1,2,3, ...\}$, endowed with a birth and death generator $\bar L$,
 namely of the form
 \bq
 \fo x\not=y\in \bar S,\qquad
 \bar L(x,y)&=&\lt\{\begin{array}{ll}
b_x&\hbox{, if $y=x+1$}\\
d_x&\hbox{, if $y=x-1$}\\
-d_x-b_x&\hbox{, if $y=x$}\\
0&\hbox{, otherwise,}\\
 \end{array}\rt.
 \eq
where $(b_x)_{x\in\ZZ_+}$ and $(d_x)_{x\in\NN}$ are the positive birth and death rates, except that $b_0=0$: 0 is the absorbing state
and the restriction of $\bar L$ to $\NN$ is irreducible.
\par
The boundary point $\iy$ is said to be an entrance boundary for $\bar L$ (cf.\ for instance Section 8.1 of the book \cite{MR1118840} of Anderson) if the following conditions are met:
\bqn{R}
 \sum_{x=1}^\iy\frac1{\pi_xb_x}\sum_{y=1}^x\pi_y&=&+\iy\\
\label{S} \sum_{x=1}^\iy\frac1{\pi_xb_x}\sum_{y=x+1}^{\iy}\pi_y&<&+\iy,\eqn
where 
\bqn{pi}
\fo x\in\NN,\qquad 
\pi_x&\df& \lt\{\begin{array}{ll}
1&\hbox{, if $x=1$}\\
\frac{b_1b_2\cdots b_{x-1}}{d_2d_3\cdots d_{x}}&\hbox{, if $x\geq 2$.}\end{array}\rt.
 \eqn
The meaning of \eqref{R} is that it is not possible (a.s.) for the underlying process $X^x$, for $x\in \ZZ_+$ to explode to $\iy$ in finite time,
while \eqref{S} says it can come back  in finite time from as close as wanted to $\iy$.\par
One consequence of \eqref{S} is that $Z\df\sum_{x\in\NN} \pi_x<+\iy$, so we can consider the probability
\bq
\fo x\in\NN,\qquad \eta(x)&\df& Z^{-1} \pi_x.\eq
Denote by $\cF$ the space of functions which vanish outside a finite subset of points from $\NN$
and by $K$ the restriction of the operator $\bar L$ to $\cF$.
It is immediate to check that $K$ is symmetric on $\LL^2(\eta)$.
Thus we can  consider its Freidrich's extension (see e.g.\ the book of Akhiezer and Glazman \cite{MR615737}), still denoted $K$, which is a self-adjoint operator in $\LL^2(\eta)$.
The fact that $\iy$ is an entrance boundary ensures indeed that such a self-adjoint extension is unique.
It is furthermore known that the spectrum of $-K$ only consists of eigenvalues of multiplicity one, say the $(\lambda_n)_{n\in\ZZ_+}$
in increasing order, see for instance Gong, Mao and Zhang \cite{MR2993011}. Let $\varphi$ be an eigenvector associated to the eigenvalue $\lambda_0>0$ of $-K$.
As in \eqref{wila}, since the absorbing point is only reachable from 1, we also introduce
\bq \lambda'_0&\df& \lambda_0(\NN\setminus\{1\}),\eq
which is the first eigenvalue of the restriction of $-K$ to functions which vanish at 1.
\par
We can now state the extension of Proposition \ref{pro2}:
\begin{theo}\label{theo}
Under the above assumptions,
we have $\lambda'_0>\lambda_0$ and 
\bqn{tracefinie}
\sum_{n\in\ZZ_+}\frac{1}{ \lambda_n}&<&+\iy.\eqn
In particular, we deduce that
\bq\lt( 1-\frac{\lambda_0}{\lambda'_0}\rt)\prod_{n\in \NN} \lt( 1-\frac{\lambda_0}{\lambda_n}\rt)&>&0.\eq
Up to a change of sign, the eigenvector $\varphi$ is increasing on $\ZZ_+$
(with the convention $\varphi(0)=0$).
It is furthermore bounded and  its amplitude satisfies:
 \bq
 \frac{\sup_{x\in\NN} \varphi(x)}{\inf_{y\in\NN} \varphi(y)}&=& \frac{\lim_{x\ri\iy} \varphi(x)}{ \varphi(1)}\\
 &\leq & \lt(\lt( 1-\frac{\lambda_0}{\lambda'_0}\rt)
\prod_{n\in \NN} \lt( 1-\frac{\lambda_0}{\lambda_n}\rt)\rt)^{-1}.\eq
\end{theo}
\par
There is a classical converse result, showing that the criterion of entrance boundary is in some sense optimal for effective absorption at 0 and  boundedness of $\varphi$.
It is typically based on the Lyapounov function approach of convergence of Markov processes (cf.\ the book of Meyn and Tweedie \cite{MR2509253}), Proposition \ref{pro3}  below gives a more precise statement for an example.
\par
Let $\bar L$ be a birth and death generator on $\ZZ_+$, absorbing at 0 and irreducible on $\NN$.
It is always possible to associate to it the minimal Markov processes $X^x\df(X_t)_{0\leq t<\sigma_\iy}$, starting from $x\in \ZZ_+$ and defined up to the explosion time $\sigma_\iy$.
These are constructed in the following  probabilistic way, where all the used random variables are independent (conditionally to the parameters
entering in the definition of their laws).
We take $X^x_t=x$ for $0\leq t<\sigma_1$, where $\sigma_1$ is distributed according to an exponential variable of parameter $\lve \bar L(x,x)\rve$
(if $x=0$, $\bar L(0,0)=0$, so $\sigma_1=+\iy=\sigma_\infty$, namely the trajectory stays at the absorbing point $0$).
Next, if $x\not= 0$, the position $X^x_{\sigma_1}=y$ is chosen according to the distribution $(L(x,y)/\vert L(x,x)\vert)_{y\in \bar S\setminus\{x\}}$.
The process stays at this position for $t\in[\sigma_1,\sigma_2)$, where $\sigma_2\df\sigma_1+\cE_2$, with $\cE_2$ an exponential variable of parameter $\lve \bar L(X^x_{\sigma_1},X^x_{\sigma_1})\rve$.
If $X^x_{\sigma_1}\not= 0$, the next position $X^x_{\sigma_2}=y$ is chosen according to the distribution $(L(X^x_{\sigma_1},y)/\vert L(X^x_{\sigma_1},X^x_{\sigma_1})\vert)_{y\in \bar S\setminus\{X^x_{\sigma_1}\}}$.
This procedure goes on up to the time $\sigma_\iy\df\lim_{n\ri+\iy}\sigma_n$ (by convention $\sigma_\iy=+\iy$ if one of the $\sigma_n$, $n\in\NN$, is infinite, which a.s.\ means that 0 has been reached).\par
We consider again the first hitting times $\tau^x_y$ defined in \eqref{tauxy}, now  for $x,y\in\ZZ_+$.
\begin{pro}\label{pro3}
Assume on one hand,  that there exist $x\in\NN$ such that 
 $\tau^x_0$ is a.s.\ finite, namely  the process $X^x$ a.s.\ ends up being absorbed at 0. Then this is true for all $x\in\ZZ_+$.
On the other hand,  that there exist a positive number $\lambda>0$ and a positive function $\varphi$ on $\NN$, with finite amplitude $a_\varphi<+\iy$,
which satisfy $K[\varphi]\leq -\lambda \varphi$.
Then $\iy$ is an entrance boundary for $\bar L$.
\end{pro}\par
Condition \eqref{S} (coming back from infinity in finite time) for 
birth and death processes satisfying \eqref{R} (non explosion) and admitting a positive generalized eigenvector (i.e.\ not necessarily belonging to $\LL^2(\eta)$) associated to a positive eigenvalue of $-K$ is also known to be equivalent to the uniqueness of the quasi-invariant probability distribution, see Theorem 3.2 of Van Doorn \cite{MR1133722} (or
 Theorem 5.4 of the book  \cite{MR2986807} of Collet,  Mart{\'{\i}}nez and San
              Mart{\'{\i}}n).
Thus the quantitative reduction (through the amplitude $a_\varphi$) of convergence to quasi-stationarity to the convergence to equilibrium presented in \cite{diaconis:hal-01002622} can  be applied to such
 birth and death processes, 
  if and only if  they admit  a unique quasi-invariant distribution.
  \par
  The uniqueness of the quasi-stationary probability was characterized in a general setting by Champagnat and Villemonais
\cite{2014arXiv1404.1349C}. It appears that $a_\varphi$ may be infinite in this situation, in particular if diffusion processes are considered
(then $\inf \varphi=0$). 
  \par\me
 The paper is constructed according to the following plan.
 In the next section, Proposition \ref{pro0} is recovered along with a probabilistic interpretation of the first Dirichlet eigenvector $\varphi$.
As a consequence, Propositions \ref{pro1} and \ref{pro2} are obtained in Section \ref{pasa}.
The situation of denumerable absorbing at 0 birth and death processes is treated in Section \ref{dss}, where an example is given.

\section{Probabilistic interpretation of $\varphi$}\label{piov}

Our main purpose here is to recover the stochastic representation of the ratio of the first Dirichlet eigenvector $\varphi$
given in Proposition \ref{pro0}. This is due to Jacka and Roberts \cite{MR1363332},
who deduce it from the corresponding discrete time result proven by Seneta \cite{MR719544}. 
Since these authors work with denumerable state spaces, for the sake of simplicity and completeness,
we present here a direct proof for finite state spaces.
\par\me
We start by recalling three simple and classical results. Consider $\cP(S)$ the set of probability measures
on $S$.
Generalizing \eqref{tauxy}, let us define, for any initial distribution $\mu\in\cP(S)$ and for any $y\in \bar S$,
\bq
\tau^\mu_y&\df& \inf\{t\geq 0\st X^\mu_t=y\},\eq
where $(X^\mu_t)_{t\geq 0}$ is a càdlàg Markov process generated by $\bar L$ and starting from $\mu$.
\begin{lem}\label{finifini}
For any $\lambda\geq 0$, we have 
\bq
\ex \mu\in\cP(S)\st \EE[\exp(\lambda\tau^\mu_\iy)]<+\iy&\Leftrightarrow& \fo  \mu\in\cP(S),\ \EE[\exp(\lambda\tau^\mu_\iy)]<+\iy.\eq
\end{lem}
\proof
It is sufficient to consider the direct implication, the reverse one being obvious.
Since for any $\mu\in\cP(S)$, we have
\bq
\EE[\exp(\lambda\tau^\mu_\iy)]&=&\sum_{x\in S}\mu(x)\EE[\exp(\lambda\tau^x_\iy)],\eq
we just need to check that \bq
 \fo x,y\in S,\qquad  \EE[\exp(\lambda\tau^x_\iy)]<+\iy&\Leftrightarrow&  \EE[\exp(\lambda\tau^y_\iy)]<+\iy,\eq
 namely
 \bqn{implicxy}
 \fo x,y\in S,\qquad  \EE[\exp(\lambda\tau^x_\iy)]<+\iy&\Rightarrow&  \EE[\exp(\lambda\tau^y_\iy)]<+\iy.\eqn
For given $x,y\in S$, let $\gamma=(\gamma_0,\gamma_1, ..., \gamma_l)$ be a path in $S$ going from $x$ to $y$
and satisfying $\bar L(\gamma_k,\gamma_{k+1})>0$. Such a path exists, by irreducibility of $K$.
Let $A^\gamma$ be the event that the first jump of the trajectory $X^x$ is from $x=\gamma_0$ to $\gamma_1$,
that the second jump of $X^x$ is from $\gamma_1$ to $\gamma_2$, ...,
that the $l$-th jump of $X^x$ is from $\gamma_{l-1}$ to $\gamma_l$.
By the probabilistic construction of $X^x$, we have that
\bqn{rxy}
\PP[A^\gamma]&=& \frac{\bar L(\gamma_0,\gamma_1)}{\vert L(\gamma_0,\gamma_0)\vert}
\frac{\bar L(\gamma_1,\gamma_2)}{\vert L(\gamma_1,\gamma_1)\vert}\cdots \frac{\bar L(\gamma_{l-1},\gamma_l)}{\vert L(\gamma_{l-1},\gamma_{l-1})\vert}\\
\nonumber&>&0.\eqn
Using the strong Markov property of $X^x$ at the minimum time between the time of the $l$-th jump time and the absorbing time, we get
\bq
 \EE[\exp(\lambda\tau^x_\iy)]&\geq &  \EE[\un_{A^\gamma}\exp(\lambda\tau^x_\iy)]\\
 &\geq & \EE[\un_{A^\gamma}\EE[\exp(\lambda\tau^y_\iy)]]\\
 &=&\PP[A^\gamma]\EE[\exp(\lambda\tau^y_\iy)],\eq
 which implies \eqref{implicxy}.\wwtbp\par
 Define 
 \bq
 \Lambda&\df& \{\lambda\geq 0\st  \fo \mu\in\cP(S), \  \EE[\exp(\lambda\tau^\mu_\iy)]<+\iy\}.\eq
 \begin{lem}\label{nu}
 We have
 \bq
 \Lambda&=&[0,\lambda_0).\eq
 \end{lem}
 \proof
 Consider $\nu\in\cP(S)$ the quasi-stationary distribution associated to $\bar L$, namely the left eigenvector of $K$ (extended to vanish at $\iy$) associated to the eigenvalue $-\lambda_0$.
For any $t\geq 0$, the distribution of $X^\nu_t$ is
$\exp(-\lambda_0t)\nu+(1-\exp(-\lambda_0 t))\delta_\iy$. It follows that 
\bq
\fo t\geq 0,\qquad \PP[\tau^\nu_\iy>t]&=&\PP[X^\nu_t\in S]\\
&=&\exp(-\lambda_0t),\eq
namely, $\tau^\nu_\iy$ is distributed according to the exponential law of parameter $\lambda_0$.
In particular, we have
\bq
\fo \lambda\geq 0,\qquad
\EE[\exp(\lambda\tau^\nu_\iy)]&=&\lt\{
\begin{array}{ll}
\frac{\lambda_0}{\lambda_0-\lambda}&\hbox{, if $\lambda<\lambda_0$}\\
+\iy&\hbox{, if $\lambda\geq \lambda_0$.}
\end{array}\rt.
\eq
The announced result follows from the previous lemma, showing that
\bq
 \Lambda&\df& \{\lambda\geq 0\st    \EE[\exp(\lambda\tau^\nu_\iy)]<+\iy\}.\eq\wwtbp\par
 For any $\lambda\in\Lambda$, we can consider the mapping $\varphi_\lambda$ defined on $\bar S$ by
 \bq
 \fo x\in\bar S,\qquad \varphi_\lambda(x)&\df& 
 \frac{\EE[\exp(\lambda\tau^x_\iy)]}{\EE[\exp(\lambda\tau^\nu_\iy)]},\eq
 where $\nu\in\cP(S)$ is the quasi-stationary distribution of $\bar L$, whose definition was recalled in the above proof
 (but for our purpose, $\nu$ could be replaced by any other fixed distribution of $\cP(S)$).
 \begin{pro}\label{convergence}
 As $\lambda\in\Lambda$ converges to $\lambda_0$, the mapping $\varphi_{\lambda}$ converges on $S$
 to a function $\varphi$ which is a positive eigenvector associated to the eigenvalue $\lambda_0$ of $-K$.
 \end{pro}
 \proof 
We begin by checking that for fixed $\lambda\in\Lambda$, $\varphi_\lambda$ satisfies
\bqn{varphilam}
\fo x\in S,\qquad \bar L[\varphi_\lambda](x)&=&-\lambda\varphi_\lambda(x).\eqn
To simplify the notation, define
\bqn{psi}
\fo x\in\bar S,\qquad \psi_{\lambda}(x)&\df& \EE[\exp(\lambda\tau^x_\iy)],\eqn
it is sufficient to show that $ \bar L[\psi_{\lambda}]=-\lambda\psi_{\lambda}$ on $S$.
This comes from the fact that for $x\in \bar S$, the quantity $\psi_{\lambda}(x)$ can be seen as the Feynman-Kac integral with respect to 
the Markov process $X^x$ and  the potential $\lambda \un_S$.
But maybe the shortest way to deduce it  is to use the martingale problem associated to $X^x$ (for a general reference, see the book of Ethier and Kurtz \cite{MR838085}).
More precisely, consider the mapping $f$ on $\RR_+\times \bar S$ defined by
\bq
\fo (t,y)\in\RR_+\times \bar S,\qquad f(t,y)&\df& \exp(\lambda t)\psi_{\lambda} (y).\eq
There exists a local martingale $M=(M_t)_{t\geq 0}$ such that a.s.
\bq
\fo t\geq 0,\qquad
f(t,X^x_t)&=&f(0,x)+\int_0^t\pa_sf(s,X^x_s)+\bar L[f(s,\cdot)](X^x_s)\, ds+M_t.\eq
The fact that $\lambda\in\Lambda$ implies that $M$ is an actual martingale
(namely that for all $t\geq 0$, $M_t$ is integrable).
In particular, by the stopping theorem, we get that for any $t\geq 0$, $\EE[M_{t\wedge \tau^x_\iy}]=0$, so that
\bq
\EE[f(t\wedge  \tau^x_\iy, X^x_{t\wedge  \tau^x_\iy})]&=& \psi_{\lambda}(x)+\EE\lt[
\int_0^{t\wedge  \tau^x_\iy}\pa_s f(s,X^x_s)+\bar L[f(s,\cdot)](X^x_s)\, ds
\rt].\eq
\par
But the strong Markov property applied to the stopping time $t\wedge  \tau^x_\iy$ implies that
\bq\EE[f(t\wedge  \tau^x_\iy, X^x_{t\wedge  \tau^x_\iy})]
&=&\EE[\exp(\lambda\tau^x_\iy)]\\
&=& \psi_{\lambda}(x),\eq
and we get that
\bq
\EE\lt[
\int_0^{t\wedge  \tau^x_\iy}\pa_s f(s,X^x_s)+\bar L[f(s,\cdot)](X^x_s)\, ds
\rt]
&=&0.\eq
Taking into account that for any $s\geq 0$ and $y\in \bar S$, $\pa_sf(s,y)=\lambda f(s,y)$,
we deduce that
\bq
\lambda \psi_{\lambda}(x)+\bar L[\psi_{\lambda}](x)&=&\lim_{t\ri 0_+}t^{-1}
\EE\lt[
\int_0^{t\wedge  \tau^x_\iy}\lambda f(s,X^x_s)+\bar L[f(s,\cdot)](X^x_s)\, ds
\rt]\\
&=&0,\eq
which amounts to \eqref{varphilam}.\par
Of course, the $\lambda\in \Lambda$ are not Dirichlet eigenvalues of $\bar L$, because $\varphi_\lambda(\iy)\not =0$:
\bq
\varphi_\lambda(\iy)&=&\frac{1}{\EE[\exp(\lambda\tau^\nu_\iy)]},\eq
but as $\lambda\in\Lambda$ goes to $\lambda_0$, this expression converges to zero.
Furthermore, if for $x,y\in S$, we call $r_{x,y}$ the r.h.s.\ of \eqref{rxy} and
\bq
r&=&\min_{x,y\in S} r_{x,y},\eq
then
\bqn{compact}
\fo \lambda\in\Lambda,\,\fo x,y\in S,\qquad r\ \leq \ \frac{\varphi_\lambda(y)}{\varphi_\lambda(x)}\ \leq \ r^{-1}.\eqn
Thus we can find a sequence $(l_n)_{n\in \NN}$ of elements of $\Lambda$ converging to $\lambda_0$
such that $\varphi_{l_n}$ converges toward a function $\varphi$ on $\bar S$, positive on $S$.
According to the previous observation $\varphi(\iy)=0$ and
taking the limit in \eqref{varphilam}, we get
\bq \bar L[\varphi]&=&-\lambda_0\varphi,\eq
it follows that the restriction to $S$ of $\varphi$ is a positive eigenvector associated to the eigenvalue $\lambda_0$ of $-K$.
Furthermore,
\bq \nu[\varphi]\ =\ \lim_{\lambda\ri\lambda_0-}\nu[\varphi_\lambda]\ =\ 1,\eq and this normalization entirely determines $\varphi$.
It follows that the mapping $\varphi$ does not depend on the chosen sequence $(l_n)_{n\in \NN}$. A usual compactness argument
based on \eqref{compact} shows that in fact
\bq
\lim_{\lambda\ri\lambda_0-}\varphi_{\lambda}&=&\varphi.\eq
\wwtbp
We need a last preliminary result.
\begin{lem}\label{lala}
For any $x\in S$, we have
\bq
\lambda_0(S\setminus\{x\})&>&\lambda_0,\eq
where we recall that the l.h.s.\ is the first eigenvalue of the $(S\setminus\{x\})\times (S\setminus\{x\})$
minor of $-\bar L$.
\end{lem}
Heuristically, this result says that for any fixed $x\in S$, it is asymptotically strictly easier for the underlying processes to exit $S\setminus\{x\}$ than $S$.
It is well-known in the reversible context, via the variational characterization of the eigenvalues, but we cannot use that argument here. Note also that in the trivial case where
$S$ is reduced to a singleton, by convention $\lambda_0(\emptyset)=+\iy$ and the above inequality is also true.
\proof
Fix $x\in S$ and let $\varphi^x$ be a positive eigenvector associated to the eigenvalue $-\lambda_0(S\setminus\{x\})$ of  the $(S\setminus\{x\})\times (S\setminus\{x\})$
minor of $-\bar L$. Extending $\varphi^x$ on $\bar S$ by making it vanish on $\{\iy, x\}$, we have that
\bq
\fo y\in S\setminus\{x\},\qquad \bar L[\varphi^x](y)&=&-\lambda_0(S\setminus\{x\})\varphi^x(y).\eq
Consider the set 
\bq
S'\ \df\ \{y\in S\st \varphi^x(y)=0\}\ \supset\ \{x\}.\eq
By irreducibility of $K$, there exists $x_0\in S'$ and $y_0\in S\setminus S'$ with $\bar L(x_0,y_0)>0$.
It follows that
\bq
\bar L[\varphi^x](x_0)&=&\sum_{y\in \bar S}\bar L(x_0,y)(\varphi^x(y)-\varphi^x(x_0))\\
&=&\sum_{y\in  S}\bar L(x_0,y)\varphi^x(y)\\
&\geq & \bar L(x_0,y_0)\varphi^x(y_0)\\
&>& 0.\eq
Similarly, we prove that
\bq
\fo y\in S',\qquad \bar L[\varphi^x](y)&\geq & 0\eq
(this is the maximum principle for the Markovian generator $\bar L$).\par
Let $\nu$ be the quasi-stationary measure associated to $\bar L$, already encountered in the proof of Lemma \ref{nu}.
Since $\nu \bar L=-\lambda_0\nu$, we have in particular
\bq
\nu[\bar L[\varphi^x]]&=&-\lambda_0\nu[\varphi^x].\eq
But according to the previous observations, we have
\bq
\nu[\bar L[\varphi^x]]&=&\nu[\un_{S\setminus S'}\bar L[\varphi^x]]+\nu[\un_{S'}\bar L[\varphi^x]]\\
&=&-\lambda_0(S\setminus\{x\})\nu[\un_{S\setminus S'}\varphi^x]
+\nu[\un_{S'}\bar L[\varphi^x]]\\
&=&-\lambda_0(S\setminus\{x\})\nu[\varphi^x]
+\nu[\un_{S'}\bar L[\varphi^x]]\\
&>&-\lambda_0(S\setminus\{x\})\nu[\varphi^x].\eq
It follows that 
\bq
\lambda_0(S\setminus\{x\})&>&\lambda_0.\eq
\wwtbp
We can now come to the
\prooff{Proof of Proposition \ref{pro0}}
Concerning the first equality, let us fix $x,y\in S$. We can assume that $x\not= y$,
since the equality is trivial for $x=y$. According to Proposition \ref{convergence},
it is sufficient to see that
\bqn{varphivarphi}
\lim_{\lambda\ri\lambda_0-}\frac{\varphi_\lambda(x)}{\varphi_\lambda(y)}
&=&\EE[\exp(\lambda_0\tau_y^x)\un_{\tau^x_y<\tau^x_\iy}].\eqn
Define $\tau=\tau^x_y\wedge\tau^x_\iy$. It is the exit time from $S\setminus\{y\}$ for $X^x$.
In particular, we have
\bq
\fo l\in\RR_+,\qquad \EE[\exp(l\tau)]<+\iy&\Leftrightarrow& l<\lambda_0(S\setminus\{y\}),\eq
and Lemma \ref{lala} implies that
\bqn{fini}
 \EE[\exp(\lambda_0\tau)]&<&+\iy\eqn
 For $\lambda\in\Lambda$, consider again the function $\psi_\lambda$ defined in \eqref{psi}.
Using the strong Markov property at time $\tau$, we have
\bq
\psi_\lambda(x)&=&\EE[\exp(\lambda\tau)\un_{X^x_{\tau}=y}\psi_\lambda(y)]
+\EE[\exp(\lambda\tau)\un_{X^x_{\tau}=\iy}\psi_\lambda(\iy)]\\
&=&\psi_\lambda(y)\EE[\exp(\lambda\tau)\un_{X^x_{\tau}=y}]+\EE[\exp(\lambda\tau)\un_{X^x_{\tau}=\iy}]\eq
Dividing by $\psi_\lambda(y)$, taking into account \eqref{fini} and letting $\lambda$ go to $\lambda_0$,
we get \eqref{varphivarphi}. In particular, we deduce that
\bq
a_\varphi&=&\max_{x\in S,\, y\in S}\EE[\exp(\lambda_0\tau_y^x)\un_{\tau^x_y<\tau^x_\iy}]\eq
To show the representation of $a_\varphi$ in Proposition \ref{pro0}, it is enough to check that 
$\varphi_\wedge=\min_{y\in O}\varphi(y)$.
This is a consequence of the fact that for any $y\in S$, either $y\in O$ or there exists a neighbor $z\in S$ of $y$ (namely a point satisfying $\bar L(y,z)>0$)
with $\varphi(z)<\varphi(y)$. Indeed this comes from
\bq
\sum_{z\in \bar S} \bar L(y,z)(\varphi (z)-\varphi(y))&=&-\lambda_0\varphi(y)\\
&<& 0.\eq
\wwtbp

\section{Path and spectral arguments}\label{pasa}

It will be seen here how the probabilistic representation of the amplitude $a_\varphi$ can be used to deduce
more practical estimates.\par\me
We begin with a path argument, similar in spirit  to the one already encountered in the proof of Lemma \ref{finifini}.
Let  $\gamma=(\gamma_0,\gamma_1, ..., \gamma_l)$ be a path in $S$, to which we associate  the
event $A^\gamma$ requiring that the first jump of the trajectory $X^{\gamma_0}$ is from $\gamma_0$ to $\gamma_1$,
that the second jump of $X^{\gamma_0}$ is from $\gamma_1$ to $\gamma_2$, ...,
that the $l$-th jump of $X^{\gamma_0}$ is from $\gamma_{l-1}$ to $\gamma_l$.
\begin{lem}
For any $\lambda\in[0, \min(\lve \bar L(\gamma_k,\gamma_k)\rve\st k\in\lin 0,l-1\rin))$, we have
\bq
\EE[\un_{A^\gamma}\exp(\lambda\tau^{\gamma_0}_{\gamma_l})]&=& \prod_{k\in \lin 0, l-1\rin}\frac{\bar L(\gamma_k,\gamma_{k+1})}{\lve \bar L(\gamma_k,\gamma_k)\rve-\lambda}.\eq
If $\lambda\geq \min(\lve \bar L(\gamma_k,\gamma_k)\rve\st k\in\lin 0,l-1\rin)$, the expectation in the l.h.s.\ is infinite.
\end{lem}
\proof
This result is directly based on the probabilistic construction of the trajectory $X^{\gamma_0}$.
Let us recall it: $X^{\gamma_0}$ stays at $\gamma_0$ for an exponential time of parameter $\lve \bar L(\gamma_0,\gamma_0)\rve$,
then it chooses a new position $x_1$ according to the probability $\bar L(\gamma_0,x_1)/\lve \bar L(\gamma_0,\gamma_0)\rve$.
Next it stays at $x_1$  for an exponential time of parameter $\lve \bar L(x_1,x_1)\rve$,
until it chooses a new position $x_2$ with respect to the probability $\bar L(x_1,x_2)/\lve \bar L(x_1,x_1)\rve$, etc.
To simplify the notation, denote
\bq
\fo k\in\lin 0,l-1\rin,\qquad L_k&\df& \lve \bar L(\gamma_k,\gamma_k)\rve.\eq
It follows that if $\lambda<  \min(L_k \st k\in\lin 0,l-1\rin)$,
\bq
\lefteqn{\EE[\un_{A^\gamma}\exp(\lambda\tau^{\gamma_0}_{\gamma_l})]}\\&=&
\lt(\prod_{k\in \lin 0, l-1\rin}\frac{\bar L(\gamma_k,\gamma_{k+1})}{ L_k}\rt)
\int\int\cdots\int e^{\lambda(t_0+t_1+\cdots +t_k)}\, L_0e^{-L_0t_0}dt_0\, L_1e^{-L_1t_1}dt_1\cdots  L_{l-1}e^{-L_{l-1}t_{l-1}}dt_{l-1}\\
&=&\prod_{k\in \lin 0, l-1\rin}\lt(\frac{\bar L(\gamma_k,\gamma_{k+1})}{L_k} L_k\int e^{(\lambda -L_k)t_k}\, dt_k\rt)\\
&=& \prod_{k\in \lin 0, l-1\rin}\frac{\bar L(\gamma_k,\gamma_{k+1})}{L_k-\lambda}
.\eq
The same computation shows that if for some $k\in\lin 0,l-1\rin$,
$\lambda \geq \lve \bar L(\gamma_k,\gamma_k)\rve$, then 
$\EE[\un_{A^\gamma}\exp(\lambda\tau^{\gamma_0}_{\gamma_l})]=+\iy$.
\wwtbp\par
\prooff{Proof of
Proposition \ref{pro1}}
It is now a consequence of the following observation:
if $\gamma_{y,x}$ is a path going from $y$ to $x$ in $S$, then from Proposition \ref{pro0},
we get
\bq
\frac{\varphi(y)}{\varphi(x)}&=&\EE[\exp(\lambda_0\tau_x^y)\un_{\tau^y_x<\tau^y_\iy}]\\
&\geq & \EE[\un_{A^\gamma}\exp(\lambda_0\tau^{y}_{x})]\\
&=& P(\gamma_{y,x}),\eq
according to the previous lemma 
(where the functional $P$ was defined in \eqref{Pg}).
Indeed, one would have noticed that 
\bq
\lambda_0&\leq & \min_{x\in S} \lve \bar L(x,x)\rve.\eq
Arguments similar to those given in the proof of Lemma \ref{lala} (reinterpret $\bar L(x,x)$ as the first Dirichlet eigenvalue
associated to the $\{x\}\times\{x\}$ minor of $\bar L$)
show that 
 the above inequality is strict, except if $S$ is a singleton.
In the latter case, say $S=\{x_0\}$, necessarily $y=x=x_0$ and $\gamma_{x_0,x_0}=(x_0)$, so that the product defining $P(\gamma_{x_0,x_0})$
is void, meaning that $P(\gamma_{x_0,x_0})=1$, as it should be.\par
Coming back to the general case and taking the minimum over $x\in S$ and $y\in O$, we get
\bq
a_\varphi&=&\lt(\min_{y\in O,\, x\in S} \frac{\varphi(y)}{\varphi(x)}\rt)^{-1}\\
&\leq & \lt(\min_{y\in O,\, x\in S}  P(\gamma_{y,x})\rt)^{-1},\eq
as announced.\wwtbp
One can deduce a rougher estimate, where $\lambda_0$ does not enter:
with the notation of \eqref{Pg}, define
\bq Q(\gamma)&\df& \prod_{k\in \lin 0, l-1\rin}\frac{\lve \bar L(\gamma_k,\gamma_k)\rve}{\bar L(\gamma_k,\gamma_{k+1})},
\eq
then we have
\bqn{aQ}
a_\varphi
&\leq & \max_{y\in O,\, x\in S}  Q(\gamma_{y,x}).\eqn
Let us illustrate these computations with
\begin{exas}\label{exa1}
(a) Consider an oriented finite strongly connected graph  $G$
on the vertex set $S$ and denote by $E$ the set of its oriented edges.
Let  $O$ be a non-empty subset of $S$. Let $\bar G$ be the oriented graph
on $\bar S=S\sqcup\{\iy\}$ obtained by adding to $E$ the edges $(x,\iy)$, with $x\in O$.
Let $d$ be the maximum outgoing degree of $\bar G$ and $D$ be the ``oriented diameter" of $G$.
Let $\bar L$ be the random walk generator associated to $\bar G$:
\bq
\fo x\not=y\in\bar S,\qquad \bar L(x,y)&\df& \lt\{\begin{array}{ll}
1&\hbox{, if $(x,y)\in\bar E$}\\
0&\hbox{, otherwise.}\end{array}\rt.\eq
Choosing geodesics (with respect to the ``oriented graph distance") for the underlying paths, the bound \eqref{aQ} implies that
\bq
a_\varphi&\leq & d^D.\eq
There is an easy comparison allowing weighted edges in the case above.
If the generator $\bar L$ is perturbed to another generator $\bar L$ only satisfying,
for some constants $0<r\leq R<+\iy$,
\bq
\fo x\not=y\in\bar S,\qquad \bar L(x,y)&\in& \lt\{\begin{array}{ll}
[r,R]&\hbox{, if $(x,y)\in\bar E$}\\
\{0\}&\hbox{, otherwise,}\end{array}\rt.\eq
we  end up with
\bqn{aRdrD}
a_\varphi&\leq & \lt(\frac{Rd}{r}\rt)^D.\eqn\par
(b) To see if this bound is of the right order, let us consider  a specific birth and death examples on $\bar S=\lin 0,N\rin$, with $N\in\NN$,
absorbed in 0 (namely $\iy$ in the previous notation). The only non-zero jump rates of $\bar L$ are given by
\bqn{L2}
\fo x\in \lin 1,N-2\rin,&\ & \lt\{
\begin{array}{ccc}\bar L(x,x+1)& \df& \rho\\ \bar L(x+1,x)& \df&1
\end{array}\rt.\\
\label{L2b}&&
 \hskip-40mm 
\bar L(1,0)=1,\quad \bar L(N-1,N)\ =\ \rho\ \hbox{ and }\ \bar  L(N,N-1)\ =\ 1+\rho,
\eqn
for some fixed $\rho>0$.
The  bound \eqref{aRdrD} leads to 
\bqn{exploN}
a_\varphi&\leq & \lt(\frac{2(1\vee \rho)}{1\wedge \rho}\rt)^N.\eqn
It was seen in \cite{diaconis:hal-01002622} that
\bqn{equiva}
a_\varphi
&=& \lt\{\begin{array}{ll}
 \frac{2N}{\pi}(1+\cO(N^{-2}))&\hbox{, if $\rho=1$}\\
  \frac{\rho}{\rho-1}(1+\cO(\rho^{-N}))&\hbox{, if $\rho>1$,}
\end{array}\rt.\eqn
but it can be deduced from the computations presented in Section 3.3 of \cite{diaconis:hal-01002622} that if  $\rho<1$ is fixed, 
then $a_\varphi$ explodes exponentially with respect to $N$.
So \eqref{exploN} corresponds to the right behavior of $a_\varphi$ (i.e.\ it does explode  exponentially with respect to $N$) if and only if $\rho<1$.
\end{exas}
In the previous example for $\rho\geq 1$, the path estimate does not catch the fact that either $a_\varphi$ is bounded
(for $\rho>1$) or explodes linearly with respect to $N$ (for $\rho=1$). The spectral estimates we are about to present
are more precise and we will see how to recover these behaviors of $a_\varphi$ for $\rho\geq 1$ and large $N$.
\par
\sm
We begin with a general result which can be deduced from Miclo \cite{MR2654550}. 
\begin{lem}
In the finite setting and under the reversibility assumption,
whatever $\mu\in\cP(S)$, the time $\tau^\mu_\iy$ is stochastically dominated by
the sum of independent exponential variables of respective parameters $\lambda_0$,
$\lambda_1$, ..., $\lambda_{N-1}$, where $N$ is the cardinality of $S$.
\end{lem}
\proof
Indeed, for any $k\in \lin 0, N-1\rin$,
denote by $\cL_k$ the convolution of $k+1$ exponential laws of parameters $\lambda_N$,
$\lambda_{N-1}$, ..., $\lambda_{N-k}$. It was seen in \cite{MR2654550} that the law of 
$\tau^\mu_\iy$ is a mixture of the $\cL_k$, for $k\in \lin 0, N-1\rin$,
the coefficients of the mixture depending on $\mu$ (and the coefficient of $\cL_{N-1}$
being positive).
The announced result follows from the fact that each of the laws  $\cL_k$, for $k\in \lin 0, N-2\rin$,
is clearly stochastically dominated by $\cL_{N-1}$.
\wwtbp
\par
We can now come to the 
\prooff{Proof of Proposition \ref{pro2}}
Note that
\bq
a_\varphi&=&1\vee\max_{x\in S,\, y\in O,\, y\not =x}\EE[\exp(\lambda_0\tau_y^x)\un_{\tau^x_y<\tau^x_\iy}],\eq
thus it is sufficient to show that
\bq
\max_{x\in S,\, y\in O,\, y\not= x}\EE[\exp(\lambda_0\tau_y^x)\un_{\tau^x_y<\tau^x_\iy}]&\leq &\lt(\lt( 1-\frac{\lambda_0}{\lambda'_0}\rt)
\prod_{k\in \lin N-1\rin} \lt( 1-\frac{\lambda_0}{\lambda_k}\rt)\rt)^{-1}.\eq
Fix $y\in O$, let $\wi K$ be the $(S\setminus\{y\})\times (S\setminus\{y\})$
minor of $\bar L$ and denote $\wi \lambda_0<\wi\lambda_1\leq \cdots \leq \wi\lambda_{N-2}$
the eigenvalues of $-\wi K$.
By the usual interlacing property of the eigenvalues of minors, we have
\bqn{inter}
\lambda_0< \wi \lambda_0\leq \lambda_1\leq \wi\lambda_1\leq \lambda_2 \leq \cdots \leq \lambda_{N-2}\leq \wi\lambda_{N-2}\leq \lambda_{N-1}.\eqn
The first inequality is strict, due to Lemma \ref{lala}.
According to the previous lemma, we have
that for any $x\in S\setminus\{y\}$, $\tau^{x}_{y}\wedge \tau^{x}_{\infty}$
is stochastically dominated by the sum of independent exponential variables of respective parameters $\wi\lambda_0$,
$\wi\lambda_1$, ..., $\wi\lambda_{N-2}$. It follows that
\bq
\EE[\exp(\lambda_0\tau_y^x)\un_{\tau^x_y<\tau^x_\iy}]&\leq &
\EE[\exp(\lambda_0(\tau^{x}_{y}\wedge \tau^{x}_{\infty}))]\\
&\leq & \EE[\exp(\lambda_0(\cE_0+\cdots+\cE_{N-2}))],\eq
where
$\cE_0$, ..., $\cE_{N-2}$ are independent exponential variables of respective parameters $\wi\lambda_0$,
$\wi\lambda_1$, ..., $\wi\lambda_{N-2}$.
Thus the last expectation is also equal to
\bqn{mieux}
\prod_{l\in\lin 0, N-2\rin}\EE[\exp(\lambda_0\cE_l)]
&=&\prod_{l\in\lin 0, N-2\rin}\frac{\wi\lambda_l}{\wi\lambda_l-\lambda_0}\\
\nonumber&\leq &\lt(\lt( 1-\frac{\lambda_0}{\wi\lambda_0}\rt)
\prod_{k\in \lin N-1\rin} \lt( 1-\frac{\lambda_0}{\lambda_k}\rt)\rt)^{-1}\\
\nonumber&\leq &\lt(\lt( 1-\frac{\lambda_0}{\lambda'_0}\rt)
\prod_{k\in \lin N-1\rin} \lt( 1-\frac{\lambda_0}{\lambda_k}\rt)\rt)^{-1},
\eqn
where the interlacing \eqref{inter} was used, as well as the definition of $\lambda'_0$ given in \eqref{wila}.
Proposition \ref{pro2} follows, since the above upper bound no longer depends on the choice of $y\in O$.\wwtbp\par
Let us now show how the spectral estimate can provide a better bound than the path estimate,
at least when some knowledge on the relevant eigenvalues is available.
\begin{exa}
We return to the birth and death processes presented at the end of Example \ref{exa1} with $\rho\geq 1$.
\par
We first treat the case $\rho=1$, for which 
we have seen in \cite{diaconis:hal-01002622}
that the eigenvalues of $-K$ are given by
\bq
\fo k\in \lin 0,N-1\rin,\qquad\lambda_k&=&
2(1-\cos((2k+1)\pi/(2N))).\eq
With the notation of the proof of Proposition \ref{pro2}, we have $O=\{1\}$ and the matrix $\wi K$ is the same as $K$, except that $N$ has been replaced by $N-1$.
Thus we get that
\bq
\fo k\in \lin 0,N-2\rin,\qquad\wi\lambda_k&=&
2(1-\cos((2k+1)\pi/(2(N-1)))).\eq
By using  \eqref{mieux} directly, we get the bound
\bq
a_\varphi&\leq & \prod_{l\in\lin 0, N-2\rin}\frac{\wi\lambda_l}{\wi\lambda_l-\lambda_0}\\
&=& \prod_{l\in\lin 0, N-2\rin}\lt(1-\frac{\lambda_0}{\wi\lambda_l}\rt)^{-1},\eq
and the first bound is in fact an equality, because
it is known that the time needed by a finite birth and death process to go from one boundary point to the other one is exactly
a sum of independent exponential variables whose parameters are the corresponding Dirichlet eigenvalues (see e.g.\ Fill \cite{MR2530102} or Diaconis and Miclo
\cite{MR2530103} for a probabilistic proof as well as a review of the history of this property).
In the above product, we begin by considering the first factor
\bq
1-\frac{\lambda_0}{\wi\lambda_0}&=&1-\frac{\sin^2(\pi/(4N))}{\sin^2(\pi/(4(N-1)))}
\\
&=&\frac{\sin^2(\pi/(4(N-1)))-\sin^2(\pi/(4N))}{\sin^2(\pi/(4(N-1)))}\\
&=&\frac{\sin(\pi/(4N(N-1)))\sin(\pi (2N-1)/(4N(N-1)))}{\sin^2(\pi/(4(N-1)))},
\eq
where we used the trigonometric formula
\bq
\fo a,b\in\RR,\qquad \sin^2(a)-\sin^2(b)&=&\sin(a+b)\sin(a-b).\eq
Letting $N$ go to infinity, it appears that
\bqn{equivN}
\lt(1-\frac{\lambda_0}{\wi\lambda_0}\rt)^{-1}&\sim &\frac{N}{2},\eqn
which provides us with the announced linear explosion in $N$. It remains to treat the other factors
\bq
1-\frac{\lambda_0}{\wi\lambda_k}&=&1-\frac{\sin^2(\pi/(4N))}{\sin^2((2k+1)\pi/(4(N-1)))},\eq
for $k\in\lin 1, N-2\rin$.
Taking into account that 
\bq
\lim_{N\ri\iy} \frac{\sin^2(\pi/(4N))}{\sin^2(3\pi/(4(N-1)))}&=&\frac19,\eq
and that for all $\theta\in[0,\pi/2]$, $2\theta/\pi\leq \sin(\theta)\leq \theta$, we can find a  constant
$c>0$ such that for $N$ large enough,
\bq
\fo k\in\lin 1, N-2\rin,\qquad
\frac{\sin^2(\pi/(4N))}{\sin^2((2k+1)\pi/(4(N-1)))}&\leq & \frac18\wedge \frac{c}{(2k+1)^2}.\eq
This bound and the dominated convergence theorem show
\bq
\lim_{N\ri\iy}\sum_{k\in\lin 1, N-2\rin}\ln\lt( 1-\frac{\sin^2(\pi/(4N))}{\sin^2((2k+1)\pi/(4(N-1)))}\rt)
&=&\sum_{k\in\NN}\ln\lt( 1-\frac{1}{(2k+1)^2}\rt)\ >\ -\iy.\eq
The above observations and \eqref{equiva} lead in fact to Wallis' formula:
\bq
\prod_{k\in\NN\setminus\{1\}\st k\, \mathrm{even}} 1-\frac{1}{k^2}&=&\frac{\pi}4.\eq\par
We now come to the case $\rho>1$. As remarked above, for all finite birth and death processes absorbed at 0, we have an exact formula
\bq
a_\varphi
&=& \prod_{l\in\lin 0, N-2\rin}\lt(1-\frac{\lambda_0}{\wi\lambda_l}\rt)^{-1},\eq
but to exploit it, one needs  a knowledge of the eigenvalues
$\lambda_0, \wi\lambda_0, \wi \lambda_1, \dots, \wi\lambda_{N-2}$.
The only behavior provided in \cite{diaconis:hal-01002622}, is that for large $N$
\bqn{raprho1}
\lambda_0&\sim &\frac12(\rho+1)(\rho-1)^2\frac1{\rho^{N+1}}.\eqn
Let us  show how this is sufficient to deduce that $a_\varphi$ remains bounded as $N$ go to infinity.
Indeed, we will just need an additional qualitative result about the number of nodal domains
of the corresponding eigenvectors, which is a discrete analogue of Sturm's theorem for one dimensional diffusions.\par
We begin by treating the first factor. As above, by spatial homogeneity, $\wi \lambda_0$ is just $\lambda_0$ when $N$ has been replaced by $N-1$.
It follows that 
$\wi
\lambda_0\sim \frac12(\rho+1)(\rho-1)^2\frac1{\rho^{N}}$, so that \bq
\lim_{N\ri\iy}
\lt(1-\frac{\lambda_0}{\wi\lambda_0}\rt)^{-1}&=&\lt(1-\frac{1}{\rho}\rt)^{-1},\eq
which in comparison with \eqref{equivN}, is a first indication why $a_\varphi$ should stay bounded as $N$ goes to infinity.
\par
It remains to prove that 
\bq
\limsup_{N\ri\iy}-\sum_{l\in\lin 1, N-2\rin}\ln\lt(1-\frac{\lambda_0}{\wi\lambda_l}\rt)&<&+\iy.\eq
Since we know that 
\bq
\fo l\in\lin 1, N-2\rin,\qquad \frac{\lambda_0}{\wi\lambda_l}&\leq &  \frac{\lambda_0}{\wi\lambda_0}\\
&\leq & \frac{1+\rho^{-1}}{2}\ <\ 1,\eq
for $N$ large enough,
it is sufficient to find a constant $c>0$ such that
\bq
\fo l\in\lin 1, N-2\rin,\qquad \frac{\lambda_0}{\wi\lambda_l}&\leq & c \rho^{-l},\eq
or similarly, such that
\bqn{decroissexp}
\fo l\in\lin 1, N-1\rin,\qquad \frac{\lambda_0}{\lambda_l}&\leq & c \rho^{-l}.\eqn
For given $ l\in\lin 1, N-1\rin$, let $\varphi_l$ be an eigenvector of $-\bar L$ associated to the eigenvalue $\lambda_l$ and vanishing at 0.
Since $\bar L$ is a tri-diagonal matrix,  $\varphi_l$ has $l+1$ nodal domains. 
More precisely, extend $\varphi_l$ into a continuous function on $[0,N]$ by making it affine on each of the segments $[n,n+1]$
with $n\in\lin 0, N-1\rin$. Then $\varphi_l$ has exactly $l+1$ zeros: $x_0=0<x_1<\cdots<x_l$ and it was seen in Miclo~\cite{MR2438701} that if $x_k\not\in\lin 0,N\rin$, there is a natural way
to define the jump rates $\bar L(\lceil x_k\rceil, x_k)$ and $\bar L(\lfloor x_k\rfloor, x_k)$ such we have 
\bq
\fo k\in \lin 0, l\rin,\qquad \lambda_0([x_{k},x_{k+1}])&=&\lambda_l\eq
with the convention  $x_{l+1}=N$.
Since each of the segments $[x_k,x_{k+1}]$, for $k\in \lin 0,l\rin$, contains at least one integer, we have 
$x_l\geq l$ and by consequence
\bq
\lambda_l&=& \lambda_0([x_{l},N])\\
&\geq &  \lambda_0([l,N]).\eq
Due to the spacial homogeneity of the initial generator $\bar L$, $\lambda_0([l,N])$ is the same as $\lambda_0$ where $N$ is replaced by $N-l$.
The bound \eqref{decroissexp} is now an easy consequence of \eqref{raprho1}, through the existence of a constant $C\geq 1$ (depending on $\rho>1$) such that
\bq
\fo N\in\NN,\qquad C^{-1}\rho^{-N}\ \leq \ \lambda_0\ \leq \ C\rho^{-N}.\eq
 \end{exa}

\section{Some denumerable birth and death processes}\label{dss}

This section treats  denumerable birth and death processes absorbing at 0 and with $\iy$ as entrance boundary via approximation by finite absorbing birth and death processes.
\par\me
We begin by recalling the theory of approximation of birth and death processes with $\iy$ as entrance boundary, as developed by Gong, Mao and Zhang \cite{MR2993011}.
For $N\in\NN$, consider the finite state spaces $S_N\df\lin N\rin$ and $\bar S_N\df\lin 0, N\rin$ endowed with the Markovian generator $\bar L_N$
which is the restriction of $\bar L$ to $\bar S_N$, except that $\bar L_N(N,N)=-b_{N-1}$, so that a Neumann (reflection) condition is put at $N$.
The point $0$ is still absorbing for $\bar L_N$. Denote by
\bq
\lambda_{N,0}<\lambda_{N,1}< \lambda_{N,2}< \cdots< \lambda_{N,N-1},\eq
the eigenvalues of the subMarkovian generator $K_N$, the restriction of $\bar L_N$ to $S_N$.
By convention, take $\lambda_{N,n}\df+\iy$ for $n\geq N$.
\par
Theorem 5.4 of Gong, Mao and Zhang \cite{MR2993011} asserts that for any fixed $n\in\ZZ_+$, the sequence $(\lambda_{N,n})_{N\in\NN}$
is non-increasing and that
\bqn{GMZ}
\lim_{N\ri\iy} \lambda_{N,n}&=&\lambda_n,\eqn
where $(\lambda_n)_{n\in\ZZ_+}$ are the eigenvalues of $-K$ defined in the introduction.\par
For $N\in\NN\setminus\{2\}$, let $\lambda'_{N,0}$ be the smallest eigenvalue of the restriction of $\bar L_N$ to
$\lin 2,N\rin$.
Consider the absorbing Markov generator $\bar L'$ on $\NN$,
coinciding with the restriction of $\bar L$ there, except that 1 is absorbing: $\bar L'(1,1)=\bar L'(1,2)=0$.
Applying  \eqref{GMZ} with $n=0$ and with respect to $\bar L'$, shows that
the sequence $(\lambda'_{N,0})_{N\in\NN\setminus\{1\}}$ is non-increasing and
\bqn{GMZ2}
\lim_{N\ri\iy} \lambda'_{N,0}&=&\lambda'_0.\eqn\par
These convergence properties are the main ingredient to deduce the estimate of Theorem \ref{theo} from Proposition \ref{pro2}.
We will also need the fact that the eigenvector $\varphi$ can be chosen to be positive on $\NN$ and increasing. This is well-known,
see for instance Chen \cite{MR1738551} or Miclo \cite{MR2599205}, whose arguments  can be extended to the present denumerable absorbing birth and death setting.
We present a succinct proof for the sake of completeness.
\par
For any function $f$ defined on $\NN$, consider the value
\bq
\cE(f)&=&\eta(1)d_1 f^2(1)+\sum_{x\geq 1}\eta(x)b_x(f(x+1)-f(x))^2\ \in\ \RR_+\sqcup\{+\iy\}.\eq
Then $\varphi$ is a minimizer of $\cE(f)/\eta(f^2)$ over all functions $f\in \LL^2(\eta)\setminus \{0\}$.
\par
Since $\cE(f)\leq \cE(\vert f\vert)$ for any function $f$, we can assume that $\varphi$ is non-negative, up to replacing it by $\vert \varphi\vert$.
For fixed $x\in\NN$, considering the quantity $\varphi(x)$ as a variable in the ratio $\cE(\varphi)/\eta(\varphi^2)$, it appears by minimization that
$\varphi(x)\in [\min (\varphi(x-1),\varphi(x+1)),\max (\varphi(x-1),\varphi(x+1))]$ and even $\varphi(x)\in (\min (\varphi(x-1),\varphi(x+1)),\max (\varphi(x-1),\varphi(x+1)))$
if $\varphi(x-1)\not=\varphi(x+1)$. This property implies the monotonicity of $\varphi$.
Since $\varphi(0)=0$ and $\varphi$ must be positive somewhere, it follows that $\varphi$ is non-decreasing.
Consider $x_0\df\max\{x\st \varphi(x)=0\}$. From $\varphi(x)\in (0,\varphi(x+1))$, we end up with a contradiction  if $x\not=0$.
So $x_0=0$ and $\varphi$ is positive on $\NN$. The same argument shows that if there exists $x\in\NN$ such that $\varphi(x)=\varphi(x+1)$ then $\varphi(x-1)=\varphi(x)$.
By iteration it would imply that $\varphi$ vanishes identically. Thus $\varphi$ is increasing, not only non-decreasing.\par\sm
This observation is also valid for $\bar L'$: there is an eigenvector $\varphi'$ associated to the eigenvalue $-\lambda'_0$ (of $K'$, the restriction to $\NN\setminus\{1\}$
of $\bar L'$)
which is positive and increasing on $\NN\setminus\{1\}$. Indeed, $\iy$ is also an entrance boundary for $\bar L'$, so that its spectrum
consists equally of eigenvalues of multiplicity 1, in particular $\varphi'$ exists.
As a consequence we get that $\lambda'_0>\lambda_0$:
\begin{lem}\label{lalapr}
With the above notation, 
\bq
\lambda_0'-\lambda_0&=& \eta(1)b_1\frac{\varphi'(2)\varphi(1)}{\eta[\varphi'\varphi]}\ >\ 0\eq
(where $\varphi'$ is seen as function defined on $\NN$ with the convention $\varphi'(1)=0$).
\end{lem}
\proof
The result follows from the computation of $\eta[\varphi'K[\varphi]]$: by definition of $\varphi$,
\bq
\eta[\varphi'K[\varphi]]&=&-\lambda_0\eta[\varphi'\varphi].\eq
By self-adjointness of $K$ the l.h.s.\ is equal to $\eta[\varphi K[\varphi']]$.
We remark that for $x\in\NN$,
\bq
K[\varphi'](x)&=& K'[\varphi'](x)+ b_1\varphi'(2)\un_{\{1\}}(x)\\
&=&-\lambda_0'\varphi'(x)+ b_1\varphi'(2)\un_{\{1\}}(x)\eq
(by convention, $K'[\varphi'](1)=0$), so that by multiplication by $\varphi(x)$ and  integration with respect to $\eta$,
\bq
\eta[\varphi K[\varphi']]
&=&-\lambda_0'\eta[\varphi'\varphi]+\eta(1) b_1\varphi'(2)\varphi(1).\eq
\wwtbp
\par
Let us next check the second assertion of Theorem \ref{theo}. 
\begin{lem}
Under the entrance boundary condition, \eqref{tracefinie} is valid.
\end{lem}
\proof
If we were working with an ergodic birth and death on $\ZZ_+$, this result is due to 
Mao \cite{MR2122801}.
To come back to this situation, let us consider the Markov generator $\wit L$ on $\NN$ which coincides
with $\bar L$, except that $\wit L(0,1)=1=-\wit L(0,0)$.
For this process, $\iy$ is still an entrance boundary. Let $(\wit\lambda_n)_{n\in\ZZ_+}$ be the eigenvalues of $-\wit L$. We have $\wit \lambda_0=0$ and
Mao \cite{MR2122801} has shown that
\bq
\sum_{n\in\NN}\frac{1}{\wit \lambda_n}&<&+\iy.\eq
It is well-known that the eigenvalues $(\lambda_n)_{n\in\ZZ_+}$ and $(\wit\lambda_n)_{n\in\ZZ_+}$
are interlaced:
\bq
\wit \lambda_0<\lambda_0\leq \wit\lambda_1\leq \lambda_1\leq \cdots\eq
(see for instance Miclo \cite{miclo:hal-00957019} where this kind of comparison was extensively used).
This implies the validity of \eqref{tracefinie}.\wwtbp
We can now readily end the 
\prooff{Proof of the bound of Theorem \ref{theo}}
For $N\in\NN$, let $\varphi_N$ be an eigenvector associated with the eigenvalue $\lambda_{N,0}$ of $K_N$ and normalized by $\varphi_N(1)=1$.
According to Proposition~\ref{pro2}, whose reversibility assumption is satisfied, we have
\bqn{avN}
a_{\varphi_N}&\leq &\lt(\lt( 1-\frac{\lambda_{N,0}}{\lambda_{N,0}'}\rt)
\prod_{k\in \lin N-1\rin} \lt( 1-\frac{\lambda_{N,0}}{\lambda_{N,n}}\rt)\rt)^{-1}.\eqn
Let $N_0\in\NN$ be large enough so that 
\bq
\lambda_{N_0,0}&\leq &\frac{ \lambda_0+\lambda_0'\wedge \lambda_1}{2}\quad (\,>\,\lambda_0).\eq
It follows that for $N\geq N_0$,
\bq
1-\frac{\lambda_{N,0}}{\lambda_{N,0}'}&\geq & 1-\frac{\lambda_{N,0}}{\lambda_{0}'}\\
&\geq & 1-\frac{\lambda_0+\lambda_0'}{2\lambda_{0}'}\\
&=&\frac{\lambda'_0-\lambda_0}{2\lambda_{0}'}.\eq
In a similar manner, we get that for any $n\in \NN$,
\bq
1-\frac{\lambda_{N,0}}{\lambda_{N,n}}&\geq &1-\frac{\lambda_0+\lambda_1}{2\lambda_n},\eq
so that for all $N\geq N_0$,
\bq
a_{\varphi_N}&\leq & \lt(\lt( \frac{\lambda'_0-\lambda_0}{2\lambda_{0}'}\rt)
\prod_{k\in \NN} \lt( 1-\frac{\lambda_0+\lambda_1}{2\lambda_n}\rt)\rt)^{-1},\eq
which is finite because of Lemma \ref{lalapr} and \eqref{tracefinie}.\par\sm
The functions $\varphi_N$ are also increasing on $\lin N-1\rin$. Consider them as non-decreasing mappings defined on $\NN$
by taking $\varphi_N(n)=\varphi_N(N)$ for all $n\geq N$. Due to this monotonicity property and to the above uniform bound on $a_{\varphi_N}$ over $N\geq N_0$,
we can find an increasing subsequence $(N_l)_{l\in\NN}$ and a non-decreasing and bounded function $\wi\varphi$ on $\NN$ with $\wi\varphi(1)=1$
such that $\varphi_{N_l}$ converge uniformly on $\NN$ toward $\wi\varphi$ as $l$ goes to infinity.
We are then allowed to pass to the limit in the equation $K_N[\varphi_N]=-\lambda_{N,0}\varphi_N$
to get
on $\NN$,
\bq
K[\wi\varphi]=-\lambda_0\wi \varphi.\eq
Since the function $\wi\varphi$ is bounded, it also belongs to $\LL^2(\eta)$
and so it is an eigenvector associated to the eigenvalue $-\lambda_0$ of $K$. It must thus be proportional to $\varphi$.
\par
The previous considerations also enable to pass to the limit in \eqref{avN} and this ends the proof of Theorem \ref{theo}.\wwtbp\par\me
Let us recall a classical 
\prooff{Proof of Proposition \ref{pro3}}
According to Karlin and McGregor \cite{MR0094854},
   the a.s.\ absorption of the  processes $X^x$, for $x\in\NN$, is equivalent to 
   \bq
 \sum_{x=1}^\iy\frac1{\pi_xb_x}&=&+\iy,\eq
and this  divergence clearly implies that of \eqref{R}.
 \par
Similarly to the proof of Proposition \ref{convergence}, consider next the mapping $f$ on $\RR_+\times \ZZ_+$
defined by 
\bq
\fo (t,x)\in\RR_+\times\ZZ_+,\qquad f(t,x)&\df& \exp(\lambda t)\varphi(x)\eq
(as usual, we impose $\varphi(0)=0$).
By the martingale problem solved by the law of $X^x$, for $x\in\NN$, the process $M\df(M_t)_{t\geq 0}$ given by
\bq
\fo t\geq 0,\qquad
M_t&\df& f(t,X^x_t)-f(0,x)-\int_0^t\pa_sf(s,X^x_s)+\bar L[f(s,\cdot)](X^x_s)\, ds\eq
is a local martingale and even a martingale, because for any fixed $t\geq 0$, $M_t$ is bounded,
due to the assumption $a_\varphi<+\iy$.
The stopped stochastic process $(M_{t\wedge \tau_{x,1}})_{t\geq 0}$ is also a martingale, so taking expectations at time $t\geq 0$, we get
\bq
\EE[f(t,X^x_{t\wedge \tau_{x,1}})]&=&\varphi(x)+\EE\lt[\int_0^{t\wedge \tau_{x,1}}\pa_sf(s,X^x_s)+\bar L[f(s,\cdot)](X^x_s)\, ds\rt]\\
&=&\varphi(x)+\EE\lt[\int_0^{t\wedge \tau_{x,1}}(\lambda \varphi+\bar L[\varphi])(X^x_s)\,\exp(\lambda s)ds\rt]
.\eq
By assumption, we have $\lambda \varphi+\bar L[\varphi]\leq 0$ on $\NN$, so that
\bq
\EE[\varphi(X^x_{t\wedge \tau_{x,1}})\exp(\lambda (t\wedge \tau_{x,1}))]&\leq & \varphi(x),\eq
and it follows that
\bq
\EE[\exp(\lambda (t\wedge \tau_{x,1}))]&\leq &a_\varphi.\eq
Letting $t$ go to infinity, we deduce that
\bq
\EE[\tau_{x,1}]&\leq & \frac{\EE[\exp(\lambda \tau_{x,1})]-1}{\lambda}\\
&\leq & \frac{ a_\varphi}{\lambda}.\eq
But it is well-known (see for instance Paragraph 8.1 of Anderson \cite{MR1118840}) that
\bq
\EE[\tau_{x,1}]&=&\sum_{y=1}^{x-1}\frac1{\pi_yb_y}\sum_{z=y+1}^{x}\pi_z,\eq
thus letting $x$ go to infinity we obtain 
\bq \sum_{y=1}^\iy\frac1{\pi_yb_y}\sum_{z=y+1}^{\iy}\pi_z\ \leq \ \frac{ a_\varphi}{\lambda}\ <\ +\iy,\eq
namely \eqref{S} is satisfied. This ends the proof that $\iy$ is an entrance boundary for $\bar L$.\wwtbp
\par\me
Finally, let us discuss the conditions \eqref{R} and \eqref{S}:
\begin{exa}
Consider the rates given for all $n\in\ZZ_+$ by
\bq
b_{n}&\df&\lt\{ \begin{array}{ll}
1&\hbox{, if $n\geq 1$}\\
0&\hbox{, if $n=0$}
\end{array}\rt.\\
d_n&\df&\lt\{ \begin{array}{ll}
n&\hbox{, if $n\geq 1$}\\
0&\hbox{, if $n=0$.}
\end{array}\rt.
\eq
The measure $\pi$ defined in \eqref{pi} is proportional to the restriction on $\NN$ of the Poisson distribution of parameter 1:
\bqn{pi2}
\fo n\in\NN,\qquad \pi_n&=&\frac{1}{n!}.\eqn
It is easily computed that \eqref{S} is not satisfied.
\par
To transform $\iy$ into an entrance boundary, the underlying Markov process must be accelerated near $\iy$: consider the rates given for all $n\in\ZZ_+$ by
\bq
b_{n}&\df&\lt\{ \begin{array}{ll}
\ln^2(e+n)&\hbox{, if $n\geq 1$}\\
0&\hbox{, if $n=0$}
\end{array}\rt.\\
d_n&\df&\lt\{ \begin{array}{ll}
n\ln^2(e-1+n)&\hbox{, if $n\geq 1$}\\
0&\hbox{, if $n=0$.}
\end{array}\rt.
\eq
The measure $\pi$ is not modified, still given by \eqref{pi2}. Nevertheless Conditions \eqref{R} and \eqref{S} are satisfied
and Theorem \ref{theo} can be applied.
\end{exa}
\bigskip
 \par
 \textbf{\large Acknowledgments:}\par\sm\noindent 
 The work of Persi \mbox{Diaconis} was partially supported by the CIMI (Centre International de Mathématiques et d'Informatique)
Excellence program while he was invited in Toulouse during the first semester of 2014. The work of Laurent Miclo was partially supported by
  the ANR STAB (Stabilité du comportement asymptotique d'EDP, de processus stochastiques et de leurs discrétisations).

 \bibliographystyle{plain}

\vskip2cm
\hskip70mm\box5

\end{document}